\documentclass{ws-idaqp}

\usepackage[latin1]{inputenc}

\usepackage{amsmath}
\usepackage{amssymb}
\usepackage[mathscr]{eucal}

    \usepackage{theorem}

\newenvironment{env}[2]{\begin{#1}#2\end{#1}}{}
    \newcommand{\beq}[1]{\begin{env}{equation}{#1}}
    \newcommand{\beqn}[1]{\begin{env}{equation*}{#1}}
    \newcommand{\bal}[1]{\begin{env}{align}{#1}}
    \newcommand{\baln}[1]{\begin{env}{align*}{#1}}
    \newcommand{\bga}[1]{\begin{env}{gather}{#1}}
    \newcommand{\bgan}[1]{\begin{env}{gather*}{#1}}
    \newcommand{\bflal}[1]{\begin{env}{flalign}{#1}}
    \newcommand{\bflaln}[1]{\begin{env}{flalign*}{#1}}
    \newcommand{\bmu}[1]{\begin{env}{multline}{#1}}
    \newcommand{\bmun}[1]{\begin{env}{multline*}{#1}}
    \newcommand{\bsp}[1]{\begin{env}{split}{#1}}

    \newcommand{\eeq}{\end{env}}
    \newcommand{\eeqn}{\end{env}}
    \newcommand{\eal}{\end{env}}
    \newcommand{\ealn}{\end{env}}
    \newcommand{\ega}{\end{env}}
    \newcommand{\egan}{\end{env}}
    \newcommand{\eflal}{\end{env}}
    \newcommand{\eflaln}{\end{env}}
    \newcommand{\emu}{\end{env}}
    \newcommand{\emun}{\end{env}}
    \newcommand{\esp}{\end{env}}

\newcommand{\lf}{\vspace{2ex}}

\renewcommand{\bf}[1]{\textbf{#1}}
\renewcommand{\it}[1]{\textit{#1}}
\renewcommand{\sc}[1]{\textsc{#1}}
\renewcommand{\sf}[1]{\textsf{#1}}
\renewcommand{\sl}[1]{\textsl{#1}}
\renewcommand{\tt}[1]{\texttt{#1}}
\newcommand{\hl}[1]{\bf{\it{#1}}}

\newcommand{\msf}[1]{\text{\small$\sf{#1}$}}

\newcommand{\cmc}[1]{\mathcal{#1}}
\newcommand{\eus}[1]{\mathscr{#1}}

\newcommand{\bb}[1]{\mathbb{#1}}

\newcommand{\nbd}[1]{$#1$\nobreakdash--}
\newcommand{\ol}[1]{\overline{#1}}

\newcommand{\vt}{\vartheta}

\newcommand{\family}[1]{\left(#1\right)}

\newcommand{\bfam}[1]{\bigl(#1\bigr)}
\newcommand{\Bfam}[1]{\Bigl(#1\Bigr)}
\newcommand{\AB}[1]{\langle#1\rangle}

\newcommand{\BAB}[1]{\Bigl\langle#1\Bigr\rangle}
\newcommand{\CB}[1]{\{#1\}}
\newcommand{\bCB}[1]{\bigl\{#1\bigr\}}
\newcommand{\BCB}[1]{\Bigl\{#1\Bigr\}}

\newcommand{\bSB}[1]{\bigl[#1\bigr]}

\newcommand{\LO}[1]{(#1]}

\newcommand{\set}[2][]{
    \ifthenelse{\equal{#1}{}}{
        \CB{#2}}{
        \CB{#1~|~#2}}}
\newcommand{\bset}[2][]{
    \ifthenelse{\equal{#1}{}}{
        \bCB{#2}}{
        \bCB{#1~|~#2}}}
\newcommand{\Bset}[2][]{
    \ifthenelse{\equal{#1}{}}{
        \BCB{#2}}{
        \BCB{#1~\big|~#2}}}
\newcommand{\zero}{\CB{0}}

\DeclareMathOperator{\ls}{\normalfont\msf{span}}
\DeclareMathOperator{\cls}{\ol{\ls}}
\DeclareMathOperator*{\limind}{lim\,ind}

\DeclareMathOperator{\id}{\normalfont\msf{id}}

\newcommand{\C}{\bb{C}}

\newcommand{\E}{\bb{E}}

\newcommand{\N}{\bb{N}}

\newcommand{\R}{\bb{R}}

\newcommand{\Z}{\bb{Z}}

\newcommand{\cB}{\cmc{B}}

\newcommand{\sB}{\eus{B}}

\newcommand{\sF}{\eus{F}}

\newcommand{\f}{\text{\scriptsize$\sF$}}

    \numberwithin{equation}{section}
    \renewcommand{\theequation}{\thesection.\arabic{equation}}

        \newcommand{\mnname}{Mathematical note.}
        \newcommand{\enname}{End of the note.}
        \newcommand{\definame}{Definition.}
        \newcommand{\propname}{Proposition.}
        \newcommand{\lemname}{Lemma.}
        \newcommand{\exname}{Example.}
        \newcommand{\exername}{Exercise.}
        \newcommand{\remname}{Remark.}
        \newcommand{\obname}{Observation.}
        \newcommand{\thmname}{Theorem.}
        \newcommand{\corname}{Corollary.}

        \renewcommand{\proof}{\lf\noindent\bf{Proof.~}}

    \theoremheaderfont{\normalfont\bfseries}
    \theoremstyle{change}
        \theorembodyfont{\rmfamily}
            \newtheorem{emp}{}[section]
                \newcommand{\bemp}[1][]{
                    \begin{emp}\hskip-\labelsep\bf{#1}\hskip\labelsep}
                \newcommand{\eemp}{\end{emp}}
\newtheorem{itemp}[emp]{}
                \newcommand{\bitemp}[1][]{
                    \begin{itemp}\hskip-\labelsep\bf{#1}\hskip\labelsep\normalfont\itshape}
                \newcommand{\eitemp}{\end{itemp}}
            \newtheorem{mn}[emp]{\mnname}
                \newcommand{\bmn}{\begin{mn}~\begin{quotation}\renewcommand{\baselinestretch}{1}\small\noindent\ignorespaces}
                \newcommand{\emn}{\end{quotation}\hfill\bf{\enname}\end{mn}}
            \newtheorem{ex}[emp]{\exname}
                \newcommand{\bex}{\begin{ex}}
                \newcommand{\eex}{\end{ex}}
            \newtheorem{exer}[emp]{\exername}
                \newcommand{\bexer}{\begin{exer}}
                \newcommand{\eexer}{\end{exer}}
            \newtheorem{defi}[emp]{\definame}
                \newcommand{\bdefi}{\begin{defi}}
                \newcommand{\edefi}{\end{defi}}
            \newtheorem{rem}[emp]{\remname}
                \newcommand{\brem}{\begin{rem}}
                \newcommand{\erem}{\end{rem}}
            \newtheorem{ob}[emp]{\obname}
                \newcommand{\bob}{\begin{ob}}
                \newcommand{\eob}{\end{ob}}
        \theorembodyfont{\normalfont\itshape}
            \newtheorem{thm}[emp]{\thmname}
                \newcommand{\bthm}{\begin{thm}}
                \newcommand{\ethm}{\end{thm}}
            \newtheorem{prop}[emp]{\propname}
                \newcommand{\bprop}{\begin{prop}}
                \newcommand{\eprop}{\end{prop}}
            \newtheorem{cor}[emp]{\corname}
                \newcommand{\bcor}{\begin{cor}}
                \newcommand{\ecor}{\end{cor}}
            \newtheorem{lem}[emp]{\lemname}
                \newcommand{\blem}{\begin{lem}}
                \newcommand{\elem}{\end{lem}}
\newenvironment{empn}[1]{\lf\noindent\bf{#1}\ignorespaces\hskip\labelsep}{\lf}
		\newcommand{\bempn}[1]{\begin{empn}{#1}}
		\newcommand{\eempn}{\end{empn}}
		\newcommand{\bitempn}[1]{\begin{empn}{#1}\normalfont\itshape}
		\newcommand{\eitempn}{\end{empn}}
                \newcommand{\bmnn}{\begin{empn}{\mnname}~\begin{quotation}\renewcommand{\baselinestretch}{1}\small\noindent\ignorespaces}
                \newcommand{\emnn}{\end{quotation}\hfill\bf{\enname}\end{empn}}
		\newcommand{\bexn}{\begin{empn}{\exname}}
		\newcommand{\eexn}{\end{empn}}
		\newcommand{\bexern}{\begin{empn}{\exername}}
		\newcommand{\eexern}{\end{empn}}
		\newcommand{\bdefin}{\begin{empn}{\definame}}
		\newcommand{\edefin}{\end{empn}}
		\newcommand{\bremn}{\begin{empn}{\remname}}
		\newcommand{\eremn}{\end{empn}}
		\newcommand{\bobn}{\begin{empn}{\obname}}
		\newcommand{\eobn}{\end{empn}}

\renewcommand{\msf}[1]{\mathsf{#1}}

\begin{document}

\title{A Simple Proof of the Fundamental Theorem\\about Arveson Systems\thanks{2000 AMS-Subjecrt classification: 46L55, 46L53, 60G20}}
\author{Michael Skeide}

\address{Dipartimento S.E.G.e S., Universit\`a\ degli Studi del Molise, Via de Sanctis\\86100 Campobasso, Italy\\E-mail: \tt{skeide@math.tu-cottbus.de}}
\maketitle


\begin{abstract}
With every \nbd{E_0}semigroup (acting on the algebra of of bounded operators on a separable infinite-dimensional Hilbert space) there is an associated Arveson system. One of the most important results about Arveson systems is that every Arveson system is the one associated with an \nbd{E_0}semigroup. In these notes we give a new proof of this result that is considerably simpler than the existing ones and allows for a generalization to product systems of Hilbert module (to be published elsewhere).
\end{abstract}



\section{Introduction}

{
An \hl{algebraic Arveson system} is a family $E^\otimes=\bfam{E_t}_{t\in(0,\infty)}$ of infinite-dimensional separable Hilbert spaces\footnote{Usually, we do not put restrictions on the dimension of $E_t$ and we include also $E_0=\C$ into the definition. Here we we find it convenient to stay as close as possible to [\refcite{Arv89}].} with unitaries $u_{s,t}\colon E_s\otimes E_t\rightarrow E_{s+t}$ that iterate associatively. Technically, an \hl{Arveson system} [\refcite{Arv89}, Definition 1.4] is the trivial bundle $(0,\infty)\times H_0$ ($H_0$ an infinite-dimensional separable Hilbert space) with its natural Borel structure equipped with a (jointly) measurable associative multiplication $((s,x),(t,y))\mapsto(s+t,xy)$ such that $x\otimes y\mapsto xy$ defines a unitary $H_0\otimes H_0\rightarrow H_0$. We put $E_t:=(t,H_0)$, define unitaries $u_{s,t}\colon E_s\otimes E_t\rightarrow E_{s+t}$ by setting $u_{s,t}((s,x)\otimes(t,y)):=(s+t,xy)$ and observe that $E^\otimes:=\bfam{E_t}_{t\in(0,\infty)}$ is an algebraic Arveson system.

For a section $x=\bfam{x_t}_{t\in(0,\infty)}$ $(x_t\in E_t)$ we shall denote by $x(t)$ the component of $x_t$ in $H_0$. A section $x$ in an Arveson system is \hl{measurable}, if the function $t\mapsto x(t)$ is measurable. The only property going beyond the structure of an algebraic Arveson system and the measurable structure of $(0,\infty)\times H_0$ that we need, is that for every two mesaurable sections $x,y$ the function $(s,t)\mapsto x(s)y(t)$ is measurable.\footnote{We do not know whether this property is equivalent to that $(0,\infty)\times H_0$ is an Arveson system but, clearly, an Arveson system fulfills this property.} From the mesurable structure of $(0,\infty)\times H_0$ alone it follows already that an Arveson system has a countable family $\BCB{\bfam{e^i_t}_{t\in(0,\infty)}\colon i\in\N}$ of measurable sections such that for every $t$ the family $\bCB{e^i_t\colon i\in\N}$ is an orthonormal basis for $E_t$ [\refcite{Arv89}, Proposition 1.15]. (Simply choose an orthonormal basis $\bCB{e^i\colon i\in\N}$ for $H_0$ and put $e^i_t:=(t,e^i)$.)\footnote{We note that these sections are also continuous for the trivial Banach bundle $(0,\infty)\times H_0$. This trivial observation has consequences for the generalization to product systems of Hilbert modules.}

With every proper \hl{\nbd{E_0}semigroup} $\vt=\bfam{\vt_t}_{t\in\R_+}$ (that is, a semigroup of unital endomorphisms $\vt_t$, proper for $t>0$) on $\sB(H)$ ($H$ an infinite-dimensional separable Hilbert space) that is \hl{normal} (that is, every $\vt_t$ is normal) and \hl{strongly continuous} (that is, $t\mapsto\vt_t(a)x$ is continuous for every $a\in\sB(H)$ and every $x\in H$) there is associated an Arveson system (which determines the \nbd{E_0}semigroup up to cocycle conjugacy).\footnote{If the \nbd{E_0}semigroups consists of automorphisms, then the associated $E_t$ would all be one-dimensional. Arveson excludes this case in the definition. While, usually, we tend to consider also the one-dimensional case, in these notes we find it convenient to stay with Arveson's convention.} There exist two proofs of the converse statetment: Every Arveson system arises as the Arveson system associated with an \nbd{E_0}semigroup. The first one obtained by Arveson in a series of papers [\refcite{Arv90a,Arv89a,Arv90}]. The second one, completely different, obtained by Liebscher [\refcite{Lie00p1}]. Both proofs are deep and difficult. It is our goal to furnish a new comparably simple proof of this \it{fundamental result about Arveson systems}.

There are two different ways of how to associate with a normal \nbd{E_0}semigroup $\vt$ on $\sB(H)$ an algebraic Arveson system, and if the \nbd{E_0}semigroup is strongly continuous, then in both cases the algebraic Arveson system is, in fact, an Arveson system. Therefore, in these notes we shall assume that, \hl{by convention}, all \nbd{E_0}semigroups are normal, while we say explicitly, if an \nbd{E_0}semigroup is assumed strongly continuous. We shall abbreviate $\id_{E_t}$ to $\id_t$.

The first construction is due to Arveson [\refcite{Arv89}, Section 2]. The (algebraic) Arveson system ${E^A}^\otimes=\bfam{E^A_t}_{t\in(0,\infty)}$ Arveson constructs from $\vt$ comes along with a \hl{nondegenerate representation} $\eta^\otimes=\bfam{\eta_t}_{t\in(0,\infty)}$ on $H$. That is, we have linear maps $\eta_t\colon E^A_t\rightarrow\sB(H)$ that fulfill $\eta_t(x_t)\eta_s(y_s)=\eta_{t+s}(x_ty_s)$ and $\eta_t(x_t)^*\eta_t(y_t)=\AB{x_t,y_t}\id_H$ and the nondegeneracy condition $\cls\eta_t(E^A_t)H=H$ for all $t>0$. Arveson showed existence of an \nbd{E_0}semigroup having a given $E^\otimes$ as associated Arveson system by contructing a nondegenerate representation of $E^\otimes$. Suppose we can find a family $w=\bfam{w_t}_{t\in(0,\infty)}$ of unitaries $w_t\colon E_t\otimes H\rightarrow H$ that satisfies $w_t(\id_t\otimes w_s)=w_{s+t}(u_{t,s}\otimes\id_H)$ (that is, $E_t\otimes(E_s\otimes H)=(E_t\otimes E_s)\otimes H$). Then $\eta_t(x_t)x:=w_t(x_t\otimes x)$ defines a nondegenerate representation of $E^\otimes$ and $\vt_t(a):=w_t(\id_t\otimes a)w_t^*$ an \nbd{E_0}semigroup that has $E^\otimes$ as associated Arveson system. We call the pair $(w,H)$ $(H\ne\zero)$ a \hl{right dilation} of $E^\otimes$ on $H$. (Putting $E_\infty:=H$, a right dilation extends the product on $E^\otimes$ to $E^\otimes\times\bfam{E_t}_{\LO{0,\infty}}$.) It is not difficult to show that every nondegenerate representation of $E^\otimes$ arises in the described way from a right dilation. Of course, if $\eta$ is the representation of ${E^A}^\otimes$ constructed by Arveson from an  \nbd{E_0}semigroup, then $\vt$ gives back that \nbd{E_0}semigroup.

The second construction is due to Bhat [\refcite{Bha96}]. The (algebraic) Arveson system ${E^B}^\otimes=\bfam{E^B_t}_{t\in(0,\infty)}$ Bhat constructs from $\vt$ comes along with a family $v=\bfam{v_t}_{t\in(0,\infty)}$ of unitaries $v_t\colon H\otimes E^B_t\rightarrow H$ that satisfies $v_t(v_s\otimes\id_t)=v_{s+t}(\id_H\otimes u_{s,t})$ (that is, $(H\otimes E^B_s)\otimes E^B_t=H\otimes(E^B_s\otimes E^B_t)$) so that $v_t(a\otimes\id_t)v_t^*$ defines an \nbd{E_0}semigroup (giving back $\vt_t(a)$). In general, if $E^\otimes$ is an (algebraic) Arveson system, we call a pair $(v,H)$ $(H\ne\zero)$ with a family $v$ of unitaries $v_t\colon H\otimes E_t\rightarrow H$ that satisfies the associativity condition a \hl{left dilation} of $E^\otimes$ on $H$. (Putting $E_\infty:=H$, a left dilation extends the product on $E^\otimes$ to $\bfam{E_t}_{\LO{0,\infty}}\times E^\otimes$.) Of course, the \nbd{E_0}semigroup $\vt$ defined by setting $\vt_t(a):=v_t(a\otimes\id_t)v_t^*$ has $E^\otimes$ as its associated Bhat system.

For our purposes it is indispensable to note that the Arveson system and the Bhat system of an \nbd{E_0}semigroup are not isomorphic but canonically anti-isomorphic (that is, they are equal as bundles, but the product of one is the opposite of the product of the other). As Tsirelson [\refcite{Tsi00p1}] has noted, they need not be isomorphic. So constructing a left dilation of an Arveson system $E^\otimes$ means producing an \nbd{E_0}semigroup that has $E^\otimes$ as associated Bhat system, while constructing a right dilation of an Arveson system $E^\otimes$ means producing an \nbd{E_0}semigroup that has $E^\otimes$ as associated Arveson system. Here our scope is to show that an Arveson system $E^\otimes$ can be obtained as the Bhat system of a strongly continuous \nbd{E_0}semigroup, that is, we wish to contruct a left dilation $(v,K)$ of $E^\otimes$ that has certain continuity properties. (By switching to the opposite of $E^\otimes$ this shows also that $E^\otimes$ may by obtained as the Arveson system associated with a strongly continuous \nbd{E_0}semigroup.) Anyway, for the proof of that the \nbd{E_0}semigroup we construct is strongly continuous we will construct also a right dilation $(w,L)$ of $E^\otimes$. In fact, a left dilation $(v,K)$ and a right dilation $(w,L)$ can be put together to obtain a unitary semigroup $u=\bfam{u_t}_{t\in(0,\infty)}$ on $K\otimes L$ by setting
\beq{\label{udef}
u_t
~:=~
(v_t\otimes\id_L)(\id_K\otimes w_t^*).
}\eeq
(Identifying $K=K\otimes E_t$ by $v_t$ and $L=E_t\otimes L$ by $w_t$, this is nothing but the ``rebracketting'' $k\otimes(x_t\otimes\ell)=(k\otimes x_t)\otimes\ell$, and illustrates that it is not always safe to use these identifications too naively; see also Skeide [\refcite{Ske05p1}] that discusses the case of spatial product systems.) Then the automorphism semigroup $\alpha=\bfam{\alpha_t}_{t\in(0,\infty)}$ defined as $\alpha_t=u_t\bullet u_t^*$ on $\sB^a(K\otimes L)$ restricts to the \nbd{E_0}semigroup $\vt$ on $\sB(K)\cong\sB(K)\otimes\id_L$. Showing that $u_t$ is strongly continuous will also show that $\vt$ is strongly (actually \nbd{\sigma}strongly) continuous. (Needless to say that if we extend $\alpha$ to all of $\R$, then $\alpha_{-t}$ $(t\in\R_+)$ defines an \nbd{E_0}semigroup on $\sB(L)\cong\id_K\otimes\sB(L)$ that has $E^\otimes$ as associated Arveson system.)

\brem
The relation between Arveson system and Bhat system of an \nbd{E_0}semigroup, between right and left dilation of an Arveson system is an instance of a far reaching duality between a von Neumann correspondence over a von Neumann algebra $\cB$ and its \it{commutant} which is a von Neumann correspondence over the commutant $\cB'$ of $\cB$. (The commutant has been introduced in Skeide [\refcite{Ske03c}] and, independently, in a version for \nbd{W^*}correpondences in Muhly and Solel [\refcite{MuSo04}].) A Hilbert space is, in particular, a correspondence over the von Neumann algebra $\C$ and the commutant of $\C$ is $\C'=\C$. In this picture, the Arveson system associated with an \nbd{E_0}semigroup turns out to be the \it{commutant system} of the Bhat system of that \nbd{E_0}semigroup. As Hilbert spaces the members of the two product systems are isomorphic but the commutant functor switches the order in tensor products. We find confirmed that the Arveson system of an \nbd{E_0}semigroup is anti-isomorphic to its Bhat system. See the survey Skeide [\refcite{Ske05a}] for a more detailed discussion of this duality. In Skeide [\refcite{Ske06p1}] we will present version of these notes for Hilbert modules. (In fact, our approach here was motivated by the wish to find a method that can be generalized to product systems of Hilbert modules.) In [\refcite{Ske06p1}] it will also come out more clearly why we insist to look rather at the Bhat system of an \nbd{E_0}semigroup than its Arveson system. For Hilbert modules a left dilation of a product system still gives rise immediately to an \nbd{E_0}semigroup via amplification, while the construction of an \nbd{E_0}semigroup from a right dilation is considerably more subtle than simple amplification. We do not have enough space to discuss this here more explicitly.

In fact, one reason why we decided to discuss the case of Hilbert space separately in these notes (and not integrated into [\refcite{Ske06p1}]) is that we wish to underline the extreme shortness of the argument, which would be obstructed by the far more exhaustive discussion in [\refcite{Ske06p1}]. Another reason is that in [\refcite{Ske06p1}] we concentrate on product systems that come shipped with a (strongly) continuous product system structure from the beginning, while we only scratch problems of measurability. The point is that every Bhat system of a strongly continuous \nbd{E_0}semigroup has a (strongly) continuous product system structure (see Skeide [\refcite{Ske03b}]) and not just a measurable one, and for (strongly) continuous product systems everything works as for Hilbert modules and without any separability assumption. The Arveson system of an \nbd{E_0}semigroup, instead, comes shipped with a product system operation that is only measurable. Only after showing that every Arveson system is the Bhat system of a strongly continuous \nbd{E_0}semigroup, we know that also the Arveson system of a strongly continuous \nbd{E_0}semigroup may be equipped with a continuous structure so that also the product system operation is continuous. See also Remark \ref{Lierem}
\erem
}

\section{The idea}

The idea how to construct a left dilation of an Arveson system as such is simple and can be explained quickly. Let $E^\otimes=\bfam{E_t}_{t\in(0,\infty)}$ be an algebraic Arveson system. To obtain a left dilation we proceed in two steps.

First, we construct a left dilation of the discrete subsystem $\bfam{E_n}_{n\in{\N}}$, that is, a Hilbert space $E$ and sufficiently associative identifications $E\otimes E_n=E$. Existence of such a left (and similarly of a right) dilation is comparably trivial, because every discrete product system of Hilbert spaces has unital units. In the the far more general case of Hilbert modules this is explained and exploited in Skeide [\refcite{Ske04p}]. For the sake of being self-contained we prove existence of left and right dilations of discrete product systems in the appendix.

In order to ``lift'' a dilation of $\bfam{E_n}_{n\in{\N}}$ to a dilation of $\bfam{E_t}_{t\in(0,\infty)}$ we consider the direct integrals $\int_a^bE_\alpha\,d\alpha$ ($0\le a<b\le\infty$). Clearly, under the identification of $x_t\in E_t$ and $x(t)\in H_0$, we find $\int_a^bE_\alpha\,d\alpha=L^2(\LO{a,b},H_0)$. We put $K=E\otimes\int_0^1E_\alpha\,d\alpha=L^2(\LO{a,b},E\otimes H_0)$. Choose $t>0$ and put $n:=\CB{t}$, the unique integer such that $t-n\in\LO{0,1}$. Then the following identifications
\bal{\notag
K\otimes E_t
&
~=~
E\otimes\family{\int_0^1E_\alpha\,d\alpha}\otimes E_t
~=~
E\otimes\int_t^{1+t}E_\alpha\,d\alpha
\\[2ex]\notag
&
~=~
\family{E\otimes E_n\otimes\int_{t-n}^1E_\alpha\,d\alpha}
\oplus
\family{E\otimes E_{n+1}\otimes\int_0^{t-n}E_\alpha\,d\alpha}
\\[2ex]\label{idea}
&
~=~
\family{E\otimes\int_{t-n}^1E_\alpha\,d\alpha}
\oplus
\family{E\otimes\int_0^{t-n}E_\alpha\,d\alpha}
~=~
K
}\eal
define a unitary $F\otimes E_t=F$. In the step from the second line to the third one we have made use of the identifications $E\otimes E_n=E$ and $E\otimes E_{n+1}=E$ coming from the dilation of $\bfam{E_n}_{n\in{\N_0}}$. Existence of the dilation of the discrete subsystem means that $E$ absorbs every tensor power of $E_1$. Just that how many factors $E_1$ have to be absorbed depends on whether $\alpha+t-n$ is bigger or smaller then $1$.

The only things that remain to be done is to show that, in a precise formulation, the identifications in \eqref{idea} iterate associatively and that the obtained \nbd{E_0}semigroup on $\sB(K)$ is strongly continuous. As explaind in the introduction, for the proof of continuity we will have to construct also a right dilation. Constructing a right dilation follows simply by inverting the order of factors in all tensor products. Note that this transition from left to right dilation is much more involved for Hilbert modules; see [\refcite{Ske06p1}].

\brem\label{Lierem}
We would like to say that our idea here is inspired very much by Liebscher's treatment in [\refcite{Lie00p1}, Theorem 8]. Also there the a major part of the construction exploits the properties of the Arveson system in the segment $\LO{0,1}$ and then puts together the segments suitably to cover the whole half-line. The more important it is to underline that the constructions \hl{are} definitely different. In Liebscher's construction the possibility to embed $\sB(E_t)$ into $\sB(E_s\otimes E_t)$ as $\id_s\otimes\sB(E_t)$ play as outstanding role. But it was our scope to produce a proof that works also for Hilbert modules and amplification of operators that act on the \hl{right} factor in a tensor product of Hilbert modules is, in general, impossible.

We also would like to mention another source of inspiration, namely, Riesz' proof of \it{Stone's theorem} on the generators of unitary groups as discussed in[\refcite{RiNa82}]. Also here the decomposition of the real line (containing the spectrum of the generator) into the product of $\LO{0,1}$ (leading to a periodic part in the unitary group) and $\Z$ (taking care for unboundedness of the generator) is crucial.
\erem


\section{Associativity}\label{algSEC}

In this section we specify precisely the operations suggested by \eqref{idea} and show that they iterate associatevely.

So let $E^\otimes=\bfam{E_t}_{t\in(0,\infty)}$ be an Arveson system with the family $u_{s,t}$ of unitaries defining the product system structure. Suppose $(E,\bfam{\breve{v}_n}_{n\in\N})$ is a dilation of the discrete subsystem $\bfam{E_t}_{t\in\N}$ of $\bfam{E_t}_{t\in(0,\infty)}$. Let $f=\bfam{f_\alpha}_{\alpha\in\LO{0,1}}$ be a section in $\bfam{E\otimes E_\alpha}_{\alpha\in\LO{0,1}}$ (that means just that $f_\alpha\in E\otimes E_\alpha$ for every $\alpha$) and choose $x_t\in E_t$. The operation suggested by \eqref{idea} sends the section $f\otimes x_t=\bfam{f_\alpha\otimes x_t}_{\alpha\in\LO{0,1}}$ in $\bfam{(E\otimes E_\alpha)\otimes E_t}_{\alpha\in\LO{0,1}}$ to a section $v_t(f\otimes x_t)=\bfam{(v_t(f\otimes x_t))_\alpha}_{\alpha\in\LO{0,1}}$ in $\bfam{E\otimes E_\alpha}_{\alpha\in\LO{0,1}}$ in such a way that $f_\alpha\otimes x_t$ ends up on $(v_t(f\otimes x_t))_{\alpha+t-n}$ (with $n:=\CB{\alpha+t}$ so that $\alpha+t-n\in\LO{0,1}$) defined by setting
\beq{\label{u_tsecdef}
(v_t(f\otimes x_t))_{\alpha+t-n}
~=~
(\breve{v}_n\otimes\id_{\alpha+t-n})(\id_E\otimes u_{n,\alpha+t-n}^*)(\id_E\otimes u_{\alpha,t})(f_\alpha\otimes x_t).
}\eeq
$\alpha\mapsto\alpha+t-n$ ($n$ depending on $\alpha$ and $t$) is just the shift modulo $\Z$ on $\LO{0,1}$ and, therefore, one-to-one.

\bprop\label{assprop}
The operations $v_t$ $(t\in(0,\infty))$ on sections iterate associatively, that is,
\beq{\label{sectass}
v_t(v_s\otimes\id_t)
~=~
v_{s+t}(\id_{(E\otimes E_\alpha)_{\alpha\in\LO{0,1}}}\otimes u_{s,t}).
}\eeq
\eprop

\proof
We must check, whether left-hand side and right-hand side of \eqref{sectass} do the same to the point $f_\alpha\otimes x_s\otimes y_t$ in the section $f\otimes x_s\otimes y_t$ for every $\alpha\in\LO{0,1}$. (Of course, this also shows that that both sides end up in $E\otimes E_{\alpha+s+t-\CB{\alpha+s+t}}$.) It is useful to note the following identities
\baln{
u_{r,s+t}^*u_{r+s,t}
&\
~=~
(\id_r\otimes u_{s,t})(u_{r,s}^*\otimes\id_t)
\\
(u_{r,s}\otimes\id_t)(\id_r\otimes u_{s,t}^*)
&\
~=~
u_{r+s,t}^*u_{r,s+t}
}\ealn
which follow from associativity of the $u_{s,t}$. Note also $(\id_E\otimes u_{n,\alpha+t-n}^*)(\id_E\otimes u_{\alpha,t})=\id_E\otimes u_{n,\alpha+t-n}^*u_{\alpha,t}$.

We put $m:=\CB{\alpha+s}$ and $n:=\CB{\alpha+s-m+t}=\CB{\alpha+s+t}-m$, so that $m+n=\CB{\alpha+s+t}$, and start with the left-hand side of \eqref{sectass}. We will surpress the elementary tensor $f_\alpha\otimes x_s\otimes y_t$ to which it is applied. The left-hand side, when applied to an argument in $(E\otimes E_\alpha)\otimes E_s\otimes E_t$ reads
\bal{\label{asscalc}\notag
&
(\breve{v}_n\otimes\id_{\alpha+s-m+t-n})(\id_E\otimes u_{n,\alpha+s-m+t-n}^*)(\id_E\otimes u_{\alpha+s-m,t})
\\\notag
&
~~~~~~~~~~~~~~~~~~~~~~~~~~~~~~~~~~~~\times\bSB{(\breve{v}_m\otimes\id_{\alpha+s-m})(\id_E\otimes u_{m,\alpha+s-m}^*)(\id_E\otimes u_{\alpha,s})\otimes\id_t}
\\\notag
&
~=~
(\breve{v}_n\otimes\id_{\alpha+s+t-m-n})(\id_E\otimes u_{n,\alpha+s+t-m-n}^*u_{\alpha+s-m,t})
\\\notag
&
~~~~~~~~~~~~~~~~~~~~~~~~~~~~~~~~~~~~~~~~~~\times(\breve{v}_m\otimes\id_{\alpha+s-m}\otimes\id_t)(\id_E\otimes u_{m,\alpha+s-m}^*u_{\alpha,s}\otimes\id_t)
\\\notag
&
~=~
(\breve{v}_n\otimes\id_{\alpha+s+t-m-n})(\breve{v}_m\otimes\id_n\otimes\id_{\alpha+s+t-m-n})
\\
&
~~~~~~~~~~~~~~~~~~~~~\times(\id_E\otimes\id_m\otimes u_{n,\alpha+s+t-m-n}^*u_{\alpha+s-m,t})(\id_E\otimes u_{m,\alpha+s-m}^*u_{\alpha,s}\otimes\id_t)
}\eal
where for exchanging the two factors in the middle we used the identity $(\id_{K_2}\otimes a)(a'\otimes\id_{L_1})=(a'\otimes a)=(a'\otimes\id_{L_2})(\id_{K_1}\otimes a)\in\sB^a(K_1\otimes L_1,K_2\otimes L_2)$ that holds for every $a\in\sB(L_1,L_2)$ and $a'\in\sB(K_1,K_2)$. In the last two factors in the last line \eqref{asscalc}, ignoring $\id_E$, we obtain
\bmun{
(\id_m\otimes u_{n,\alpha+s+t-m-n}^*u_{\alpha+s-m,t})(u_{m,\alpha+s-m}^*u_{\alpha,s}\otimes\id_t)
\\
~=~
(\id_m\otimes u_{n,\alpha+s+t-m-n}^*)(\id_m\otimes u_{\alpha+s-m,t})(u_{m,\alpha+s-m}^*\otimes\id_t)(u_{\alpha,s}\otimes\id_t)
\\
~=~
(\id_m\otimes u_{n,\alpha+s+t-m-n}^*)u_{m,\alpha+s+t-m}^*u_{\alpha+s,t}(u_{\alpha,s}\otimes\id_t)
\\
~=~
(u_{n,m}^*\otimes\id_{\alpha+s+t-m-n})u_{m+n,\alpha+s+t-m-n}^*u_{\alpha,s+t}(\id_\alpha\otimes u_{s,t}).
}\emun
Using associativity for the first two factors in that last line of \eqref{asscalc}, we find
\bmun{
(\breve{v}_n\otimes\id_{\alpha+s+t-m-n})(\breve{v}_m\otimes\id_n\otimes\id_{\alpha+s+t-m-n})
\\
~=~
(\breve{v}_{m+n}\otimes\id_{m+n}\otimes\id_{\alpha+s+t-m-n})(\id_E\otimes u_{m,n}\otimes\id_{\alpha+s+t-m-n}).
}\emun
Putting everything together, the factors containing $u_{m.n}$ and $u_{m,n}^*$ cancel out and we obtain
\baln{
&
(\breve{v}_{m+n}\otimes\id_{m+n}\otimes\id_{\alpha+s+t-m-n})(\id_E\otimes u_{m,n}\otimes\id_{\alpha+s+t-m-n})
\\
&
~~~~~~~~~~~~~~~~~~~~~~~~~~~~~~\times(\id_E\otimes u_{n,m}^*\otimes\id_{\alpha+s+t-m-n})(\id_E\otimes u_{m+n,\alpha+s+t-m-n}^*)
\\
&
~~~~~~~~~~~~~~~~~~~~~~~~~~~~~~~~~~~~~~~~~~~~~~~~~~~~~~\times(\id_E\otimes u_{\alpha,s+t})(\id_E\otimes\id_\alpha\otimes u_{s,t})
\\
&
~=~
(\breve{v}_{m+n}\otimes\id_{m+n}\otimes\id_{\alpha+s+t-m-n})(\id_E\otimes u_{m+n,\alpha+s+t-m-n}^*)
\\
&
~~~~~~~~~~~~~~~~~~~~~~~~~~~~~~~~~~~~~~~~~~~~~~~~~~~~~~\times(\id_E\otimes u_{\alpha,s+t})(\id_E\otimes\id_\alpha\otimes u_{s,t}).
}\emun
As $m+n=\CB{\alpha+s+t}$, this is exactly how $u_{s+t}(\id_{(E\otimes E_\alpha)_{\alpha\in\LO{0,1}}}\otimes u_{s,t})$ acts on $f_\alpha\otimes x_s\otimes y_t$.~\qed

\lf
The mapping $f\otimes x_t\mapsto v_t(f\otimes x_t)$ is fibre-wise unitary. Measurability of the product system structure implies that $v_t$ sends measurable sections to measurable sections. So, by translational invariance of the Lebesgue measure, $v_t$ sends square integrable sections isometrically to square integrable sections and, therefore, defines a unitary $v_t\colon F\otimes E_t\rightarrow F$. We summarize:

\bthm\label{ldthm}
$(v,K)$ is a left dilation of $E^\otimes$.
\ethm

Let $\f_{K,L}\colon K\otimes L\rightarrow L\otimes K$ denote the canonical \it{flip}  of the factors in a tensor product of two Hilbert spaces $K$ and $L$. Omitting the obvious proof, we add:

\bthm\label{rdthm}
If $(\tilde{v},L)$ is the left dilation constructed as before for the opposite Arveson system of $E^\otimes$, then $(w,L)$ defined by setting $w_t=\tilde{v}_t\circ\f_{E_t,L}$ is a right dilation of $E^\otimes$.
\ethm

\brem\label{switchrem}
Once more, we emphasize that an operation like the flip $\f$ is not available for Hilbert modules. In fact, in the module case $K$ will be a right Hilbert module, while $L$ will be a Hilbert space with a nondegenerate representation. Also, in the formulation of Theorem \ref{rdthm} there would not occur the opposite system of $E^\otimes$ but its commutant system. But, the commutant destroys continuity properties. Therefore, in a formulation for modules Theorem \ref{rdthm} must be reproved from scratch starting with a right dilation of the discrete subsystem (guaranteed in [\refcite{Ske04p}]) inverting in the preceeding construction the orders of the factors in all tensor products.
\erem

\brem
Recall that Proposition \ref{assprop} is a statement about an operation acting pointwise on sections and not just a statement almost everywhere. Therefore, if we replace the translation invariant Lebesgue measure by the translation invariant counting measure, so that the direct integrals become direct sums, then measurability of sections does no longer play any role. Therefore, also algebraic Arveson systems admit left and right dilations. Just that the dilation spaces might be nonseparable and the correponding \nbd{E_0}semigroups noncontinuous. The same remains true for algebraic product systems of Hilbert modules as long as we can guarantee (by [\refcite{Ske04p}]) dilations of the discrete subsystem; see [\refcite{Ske06p1}].
\erem

\section{Continuity}

As indicated in the introduction we show that the unitary semigroup on $K\otimes L$, defined with the help of the left and the right dilation from Theorems \ref{ldthm} and \ref{rdthm} by \eqref{udef}, is strongly continuous. This shows also that the \nbd{E_0}semigroup determined by $(v,K)$ is strongly continuous.

First of all, note that $E$ (carrying the dilation of the discrete subsystem) and, therefore, also $K$ and $L$ are separable. On a separable Hilbert space $H$ measurability and weak measurability are equivalent. For checking weak measurability of $t\mapsto x(t)\in H$ it is sufficient to check measurability of $t\mapsto\AB{y,x(t)}$ for $y$ from a total subset of $H$. Also, for checking strong continuity of a unitary semigroup it is sufficient to check weak measurability.

\bprop
$t\mapsto u_t$ is weakly measurable and, therefore, strongly continuous.
\eprop

\proof
Let $\BCB{\bfam{e^i_t}_{t\in(0,\infty)}\colon i\in\N}$ a mesaurable orthonormal basis for $E^\otimes$ (see the introduction). So
\beqn{
u_t
~=~
(v_t\otimes\id_L)\Bfam{\sum_{i\in\N}\id_K\otimes e^i_t{e^i_t}^*\otimes\id_L}(\id_K\otimes w_t^*)
~=~
\sum_{i\in\N}(v_t\otimes e^i_t)\otimes({e^i_t}^*\otimes w_t^*),
}\eeqn
where for every $x\in E_t$ we define the operator $v_t\otimes x\colon k\mapsto v_t(k\otimes x)$ in $\sB(K)$ and the operator $x^*\otimes w_t^*\in\sB(L)$ as the adjoint of $x\otimes w_t$ defined in a way analogous to the definition of $v_t\otimes x$. For every measurable section $y=\bfam{y_\alpha}_{\alpha\in\LO{0,1}}$ $(y_\alpha\in E_\alpha)$ and $x\in E$ the function $(\alpha,t)\mapsto v_t(x\otimes y\otimes e^i_t)_\alpha\in E\otimes E_\alpha$ is measurable and, clearly, square-integrable over $\LO{0,1}\times C$ for every compact interval $C$. Calculating an inner product with an element $k\in K$ means integrating the inner product $\AB{k_\alpha,v_t(x\otimes y\otimes e^i_t)_\alpha}$ over $\alpha\in\LO{0,1}$. By \it{Cauchy-Schwarz inequality} and \it{Fubini's theorem} the resulting function of $t$ is measurable (and square-integrable over every compact interval $C$). In orther words, the \nbd{\sB(K)}valued functions $t\mapsto v_t\otimes e^i_t$ are all weakly measurable. Similarly, the \nbd{\sB(L)}valued functions $t\mapsto {e^i_t}^*\otimes w_t^*$ are weakly measurable. By an application of the \it{dominated convergence theorem} we find that
\bmun{
t
~\longmapsto~
\BAB{(k\otimes\ell)\,,\,\Bfam{\sum_{i\in\N}(v_t\otimes e^i_t)\otimes({e^i_t}^*\otimes w_t^*)}(k'\otimes\ell')}
\\
~=~
\sum_{i\in\N}\AB{k,(v_t\otimes e^i_t)k'}\AB{\ell,({e^i_t}^*\otimes w_t^*)\ell'}
}\emun
is measurable for all $k,k\in K;\ell,\ell'\in L$. In conclusion $t\mapsto u_t$ is weakly measurable.~\qed

\lf
This concludes the proof of Arveson's theorem:

\bthm
For every Arveson system $E^\otimes$ there exists a strongly continuous \nbd{E_0}semigroup having $E^\otimes$ as associated Bhat system. By passing to the opposite Arveson system there exists also a strongly continuous \nbd{E_0}semigroup having $E^\otimes$ as associated Arveson system.
\ethm

\section*{Appendix}
\setcounter{section}{0}
\stepcounter{section}
\renewcommand{\theemp}{A.\arabic{emp}}

\bthm
Every discrete product system $\bfam{E_n}_{n\in\N}$ of (infinite-dimendional separable) Hil\-bert spaces $E_n$ admits a left and a right dilation on an (infinite-dimensional separable) Hilbert space.
\ethm

\proof
It is known since the defintion of product systems in [\refcite{Arv89}] that existence of a (unital) unit allows to construct easily a representation (that is, a right dilation). This does not depend on the index set $\N$ or $\R_+\backslash\zero$. We repeat here the construction of a left dilation from [\refcite{BhSk00}] reduced to the case of Hilbert spaces. It is noteworthy that for Hilbert modules existence of a unit vector constitutes a serious problem, while every nonzero Hilbert space has unit vectors in abundance.

Choose a unit vector $\xi_1\in E_1$ and define $\xi_n:=\xi^{\otimes n}\in E_n=E_1^{\otimes n}$. (The $\xi_n$ form a unit in the sense that $\xi_m\otimes\xi_n=\xi_{m+n}$.) The mappings $\xi_m\otimes\id_n\colon x_n\mapsto\xi_m\otimes x_n$ form an inductive system of isometric embeddings $E_n\rightarrow E_m\otimes E_n=E_{m+n}$. The inductive limit $E:=\ol{\limind_nE_n}$ comes shipped with unitaries $v_n\rightarrow E\otimes E_n\rightarrow E$ (namely, the limits of $u_{m,n}$ for $m\to\infty$) that form a left dilation.

For a right dilation either proceed as in Theorem \ref{rdthm} (just now for discrete index set) or repeat the preceeding construction for the isometric embeddings $\id_n\otimes\xi_m$. (Note that none of these suggestions works for Hilbert modules; cf.\ Remark \ref{switchrem}.)~\qed

\lf\noindent
\bf{Acknowledgements.}~
We would like to thank W.\ Arveson and V.\ Liebscher for intrigiung discussions. This work is supported by research funds of University of Molise and Italian MIUR (PRIN 2005).

\lf\noindent
\bf{Note added after acceptance.}~
After acceptance of this note Arveson [\refcite{Arv06}] has provided yet another short proof. In [\refcite{Ske06p4}] we show that the two proofs, ours here and Arveson's in [\refcite{Arv06}], actually lead to unitarily equivalent constructions. The discussions in [\refcite{Ske06p4}] unifies the advantages of the proof here (unitality of the endomorphisms) and in [\refcite{Arv06}] (no problems with associativity). In Skeide [\refcite{Ske06p5}] we applied Arveson's idea to the case of continuous product systems.

\newcommand{\Swap}[2]{#2#1}\newcommand{\Sort}[1]{}
\providecommand{\bysame}{\leavevmode\hbox to3em{\hrulefill}\thinspace}
\providecommand{\MR}{\relax\ifhmode\unskip\space\fi MR }
\providecommand{\MRhref}[2]{%
  \href{http://www.ams.org/mathscinet-getitem?mr=#1}{#2}
}
\providecommand{\href}[2]{#2}




\end{document}